\documentclass[12pt,letterpaper]{article}
\usepackage[letterpaper,margin=1in]{geometry}

\usepackage{amsmath}
\usepackage{amsfonts}
\usepackage{amsthm}
\usepackage{amstext}
\usepackage{amssymb}
\usepackage{txfonts}
\usepackage{enumerate}
\usepackage{graphicx}
\usepackage{setspace}
\usepackage{epsfig}
\usepackage{color}
\usepackage{float}
\usepackage[hang,small]{caption}

\begin{document}

\title{There are Thin Minimizers of the $L^1TV$ Functional}
\author{Benjamin Van Dyke \footnote{Department of Mathematics, Washington State University} \and Kevin R. Vixie \footnotemark[\value{footnote}]}
\date{}
\maketitle

\begin{abstract}
In this paper we show the surprising results
that while the local reach of the boundary of an $L^1TV$
minimizer is bounded below by $\frac{1}{\lambda}$, the
global reach can be smaller. We do this by demonstrating
several example minimizing sets not equal to the union of
the $\frac{1}{\lambda}$-balls they contain.
\end{abstract}


\section{Introduction}

The $L^1TV$ functional introduced and studied in \cite{chan2005} is defined to be
\begin{equation}
\label{l1tv}
F(u)=\int |\nabla u| dx + \lambda \int |u-f| dx,
\end{equation} 
where $f$ and $u$ are functions from $\mathbb{R}^n$ to $\mathbb{R}$.  If the input function $f$ is binary, Chan and Esedoglu observed that the functional reduces to:
\begin{equation}
\label{CE}
E(\Sigma,\lambda)=\mbox{Per}(\Sigma)+\lambda|\Sigma\triangle\Omega|
\end{equation}
where $\Sigma$ is the support of the function $u=\chi_{\Sigma}$, $\mbox{Per}(\Sigma)$ is the perimeter of the set $\Sigma$, $\triangle$ denotes the symmetric difference, and $\Omega$ is the support of the binary data $f=\chi_{\Omega}$.   In this paper, we will give examples to show that when $f$ is the characteristic function of a set $\Omega$, the minimizer $\Sigma_*$ of the functional is sometimes the set $\Omega$ itself, instead of the union of all the $\frac{1}{\lambda}$-balls it contains, even though there are parts of $\Omega$ that cannot contain such a ball.

We present these examples not only because they have interesting properties, but also to illustrate useful applications for many of the results found about minimizers of (\ref{l1tv}), specifically those found in \cite{Allard1,vixie}.  Furthermore, one can use these examples to test computational schemes for minimizing (\ref{l1tv}).

In \cite{Allard1,Allard2,Allard3}, Allard used techniques from geometric measure theory to produce a study of minimizers for a class of functionals that include (\ref{l1tv}).  In \cite{Allard2}, he shows that if $\Omega$ is convex, then the minimizer of (\ref{l1tv}) is either the empty set or the union of all $\frac{1}{\lambda}$-balls contained in $\Omega$.  This leads us to question:  Under what circumstances can the condition of convexity be relaxed in order to obtain the same result?  The examples in this paper explore this idea and show that the answer may be difficult.  To construct these examples we will rely on many of the results from \cite{Allard1}.  These results, as they relate to minimizers $\chi_{\Sigma_*}$ of (\ref{l1tv}) with $\Omega\in\mathbb{R}^2$, are summarized below (for the full results see \cite{Allard1}).

\begin{figure}[h!]
\centering
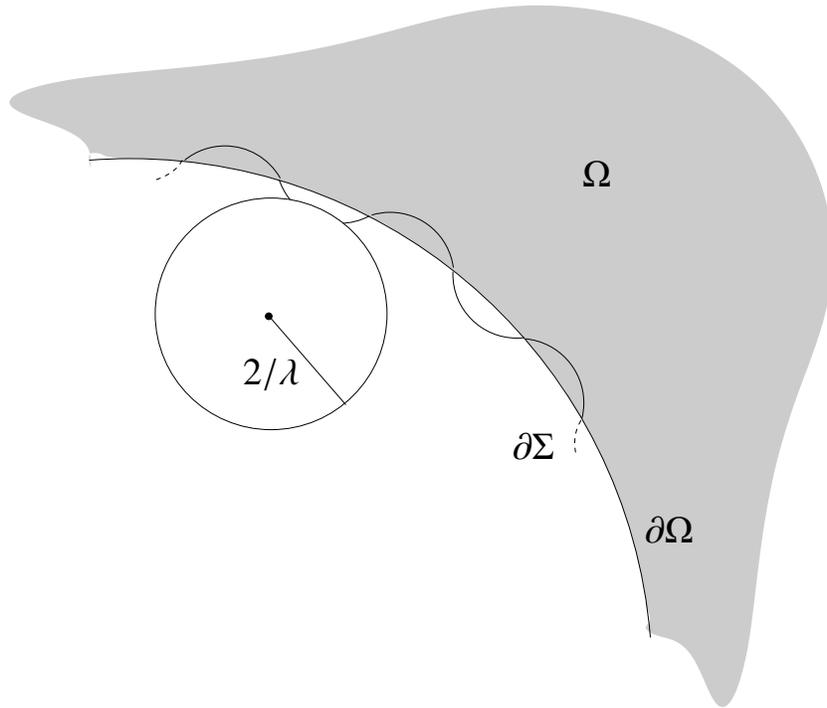
\caption{Example of the boundary of a possible minimizer $\Sigma$ that does not meet $\Omega$ tangentially and illustration of a ball of radius $2/\lambda\subset\Omega^c$ that can be used to find a new minimizer, as described in \cite{vixie}.}\label{badsigma2}
\end{figure}

\begin{enumerate}[$\bullet$]
\item The boundary $\partial\Sigma_*$ is of class $C^{1,1}$.
\item The curvature of $\partial\Sigma_*$ is bounded above by $\lambda$.
\item $\partial\Sigma_*$ differs from $\partial\Omega$ in arcs of $\lambda$-curvature.
\item These arcs subtend angles of not more than $\pi$ radians.
\item $\Sigma_*$ is contained in the closed convex hull of $\Omega$.
\end{enumerate}
Note:  If $\Sigma_*\neq\emptyset$ then $\Sigma_*$ and $\Omega$ must share part of their boundaries, otherwise the third and fourth results would be violated.  These results also imply that $\partial\Sigma_*$ and $\partial\Omega$ must meet tangentially or $\partial\Sigma_*$ is comprised of arcs of $\lambda$-curvature that meet tangentially to one another at points on $\partial\Omega$ (see figure \ref{badsigma2}).

The following theorem from \cite{vixie} can be used to eliminate the second case, illustrated in figure \ref{badsigma2}, for the three possible choices of $\Omega$ given in this paper.

\newtheorem*{theorem}{Theorem}
\begin{theorem}
Let $\Omega$ be a bounded, measurable subset of $\mathbb{R}^2$. Let $\Sigma$ be any minimizer of (\ref{CE}). Assume that a ball, $B_{2/\lambda}$, of radius $2/\lambda$ lies completely in $\Omega$. Then $B_{2/\lambda}\bigcup\Sigma$ is also a minimizer. Moreover, if $B_{2/\lambda}\subset\Omega^c$, then $(B_{2/\lambda}\bigcup\Sigma^c)^c$ is also a minimizer.
\end{theorem}

For the examples that follow, we will always be able to find a ball of radius $2/\lambda$ contained within $\Omega^c$.   Thus if a minimizer did exist with boundary as in figure \ref{badsigma2}, we would be able to find the new minimizer $(B_{2/\lambda}\bigcup\Sigma^c)^c=\Sigma\setminus B_{2/\lambda}$.  This new minimizer would then have boundary containing an arc of radius $\lambda/2$-curvature  (see figure \ref{badsigma2}) contradicting the third result listed above.  Note, it doesn't matter whether $\Sigma$ is to the left or right of $\partial\Sigma$ in figure \ref{badsigma2}, the contradiction would still be obtained.  Consequently, all the examples that follow will have minimizers with boundaries that meet $\partial\Omega$ tangentially.  We can then use the five results listed above to find the set of all possible minimizers, $\{\Sigma_i\}$,  for a given set $\Omega$ and then compute and compare their values for (\ref{CE}) to find the actual minimizer.  In each case we generate a large set of examples for which the set $\Omega$ is the minimizer despite the facts that $\bigcup\{B_{1/\lambda}(x)\subset\Omega\}\neq\Omega$ and there are parts of $\Omega$ that cannot contain such a ball.


\section{Non-concentric Annulus}

\begin{figure}[h!]
\centering
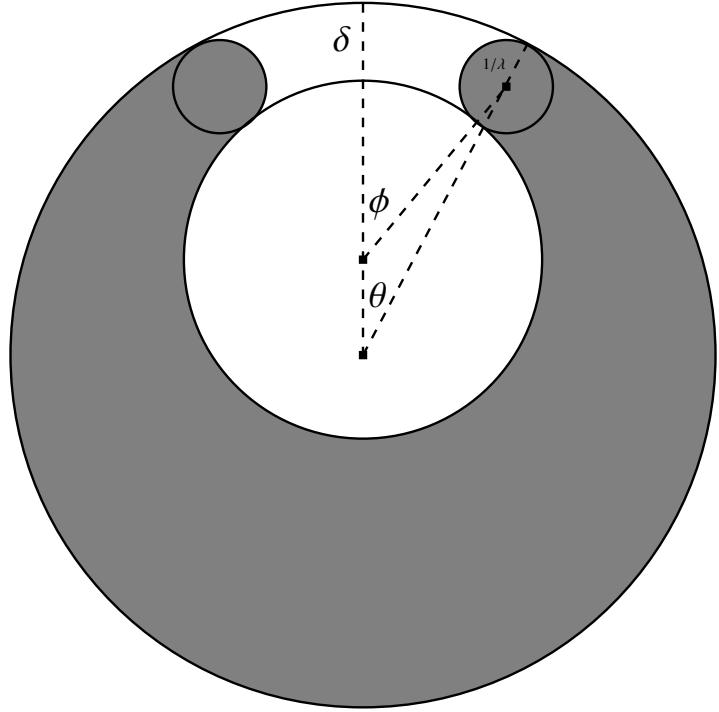
\caption{Non-concentric Annulus with angles $\phi$ and $\theta$}\label{circann}
\end{figure}

For the first example we take $\Omega$ to be the region contained between two non-concentric circles of radii $R$ and $r$ with $2/\lambda<r<R$ and minimum distance between the two circles being $\delta$ (shown in figure \ref{circann}) and compute (\ref{CE}) for the five choices of $\Sigma_i$:\\
$\Sigma_1=D_R$ (The large outer disc),\\
$\Sigma_2=\emptyset$,\\
$\Sigma_3=\Omega$,\\
$\Sigma_4=\bigcup\{B_{1/\lambda}(x)\subset\Omega\}$, and\\
$\Sigma_5=D_r$ (The small inner disc).

Since $D_R$ is equal to the closed convex hull of $\Omega$, all possible minimizers must be a subset of this disc.  This fact, along with the fact that $D_R$, $\emptyset$, $\Omega$, and $D_r$ don't violate the curvature requirement and don't have boundaries not contained in $\partial\Omega$, explain why they are considered as possible minimizers.  Since any minimizer must meet $\Omega$ tangentially, the only other possibilities would be combinations of subsets of $\Omega$ and $D_r$ that meet the annulus tangentially with arcs of $\lambda$-curvature of not more than $\pi$ radians.

We can consider both cases separately and then take unions of the possible subsets of $\Omega$ and $D_r$.  The above mentioned requirement would disallow any nonempty proper subset of $D_r$ because the restriction $\lambda>2/r$ would make it impossible for an arc of $\lambda$-curvature of not more than $\pi$ radians to meet $D_r$ tangentially.  We can now turn our attention to possible subsets of $\Omega$.  If we only consider arcs of $\lambda$-curvature of not more than $\pi$ radians that meet $\Omega$ tangentially, then we are limited to arcs that are simultaneously tangent to both circles comprising the annulus.  The only possible set that can satisfy this requirement is $\bigcup\{B_{1/\lambda}(x)\subset\Omega\}$.  Now any minimizer besides those listed previously must be either $\bigcup\{B_{1/\lambda}(x)\subset\Omega\}$ or the union of this set with one of the others.  This only leaves the possibility of $\bigcup\{B_{1/\lambda}(x)\subset\Omega\}\bigcup B_r$, but this set does not have $C^{1,1}$ boundary and so is not a possible candidate.

This leads to the following equations derived from (\ref{CE}):

\begin{equation}
\label{E1}
E(\Sigma_1,\lambda)=2\pi R+\lambda\pi r^2,\end{equation}
\begin{equation}
\label{E2}
E(\Sigma_2,\lambda)=\lambda\pi(R^2-r^2),\end{equation}
\begin{equation}
\label{E3}
E(\Sigma_3,\lambda)=2\pi(R+r),\end{equation}
\begin{equation}
\label{E4}
E(\Sigma_4,\lambda)=2(\pi-\phi)r+2(\pi-\theta)R+\frac{\pi-\phi+\theta}{\lambda}+\lambda R^2\theta-\lambda r^2\phi-(R-r-\delta)(\lambda R-1)\sin\theta\end{equation}
\begin{equation}
\label{E5}
E(\Sigma_5,\lambda)=2\pi r+\lambda\pi R^2\end{equation}
Where $\phi=\cos^{-1}\left(\frac{(R-\frac{1}{\lambda})^2-(R-r-\delta)^2-(r+\frac{1}{\lambda})^2}{2(R-r-\delta)(r+\frac{1}{\lambda})}\right)$ and $\theta=\cos^{-1}\left(\frac{(R-\frac{1}{\lambda})^2+(R-r-\delta)^2-(r+\frac{1}{\lambda})^2}{2(R-r-\delta)(R-\frac{1}{\lambda})}\right)$.  The angles $\phi$ and $\theta$ are the angles between the vertical axis and the lines from the centers of the two circular boundaries of $\Omega$ to the center of one of the two $\frac{1}{\lambda}$-balls that are tangent to both of the circular boundaries, as shown in figure \ref{circann}.

\begin{figure}[h!]
\centering
\includegraphics[width=\textwidth]{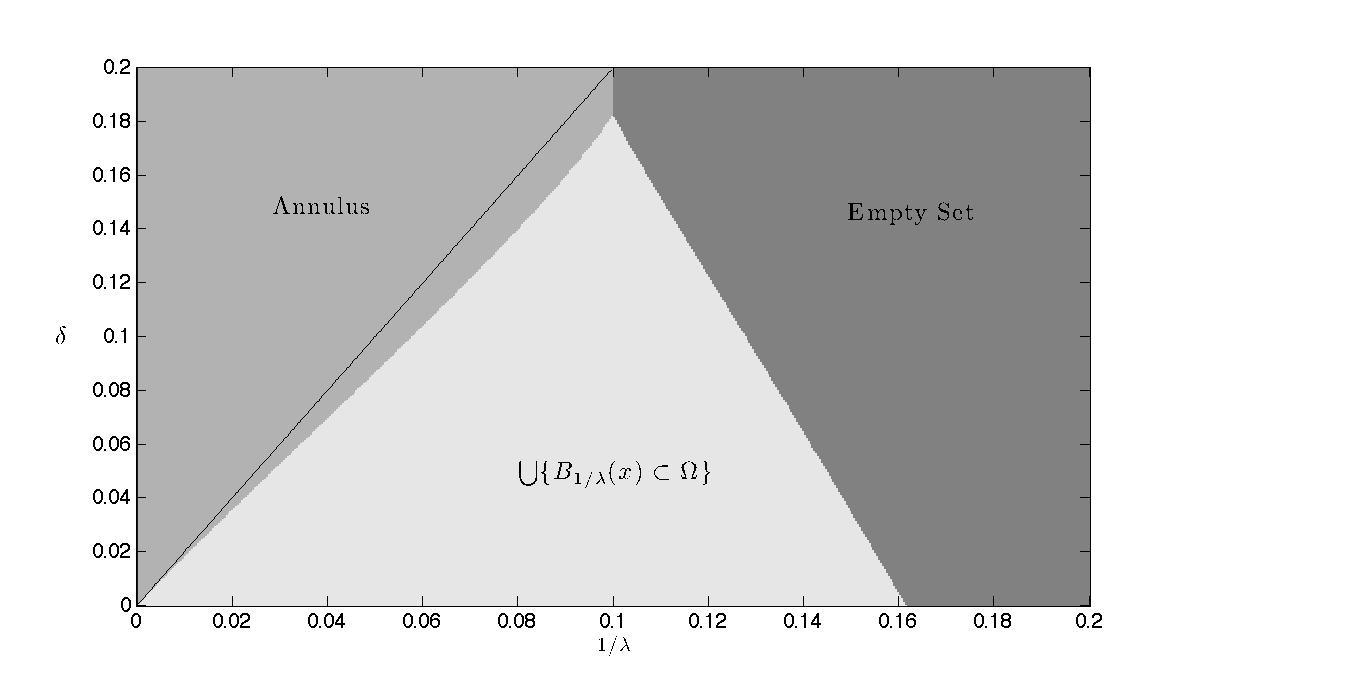}
\caption{This diagram shows the $\Sigma_i$ with minimum value in (\ref{CE}) for the non-concentric annulus.}\label{circphase}
\end{figure}

From equations (\ref{E1})-(\ref{E5}), many examples can be created.  For illustrative purposes, we have chosen $R=1$ with $r=.8$ and let $\delta$ and $\lambda$ vary.  We can then compute equations (\ref{E1})-(\ref{E5}) for any values of $\delta$ and $\lambda$ and determine which equation has minimum value, thus indicating the minimizer.  Figure \ref{circphase} summarizes the results for many such choices of $\delta$ and $\lambda$.  Since the value obtained from equation (\ref{E4}) is not meaningful for all values of $\lambda$, we have indicated with a curve on the figure where the two tangent $\frac{1}{\lambda}$-balls pictured in figure \ref{circann} would touch.  Anything to the left of this curve would indicate either the two balls pictured overlap, coincide, or $1/\lambda$ is too small for such a ball to touch both boundaries and would make (\ref{E4}) meaningless.  In the first case $\bigcup\{B_{1/\lambda}(x)\subset\Omega\}$ would not have $C^{1,1}$ boundary and so is not a possible minimizer and in the latter two cases $\bigcup\{B_{1/\lambda}(x)\subset\Omega\}=\Omega$ and so again shouldn't be considered.  Consequently, we see that to the right of the curve there is a significant region where the entire annulus obtains a lower value for (\ref{CE}) than the union of $\frac{1}{\lambda}$-balls, giving the desired examples.

\noindent Notes:  When the two tangent $\frac{1}{\lambda}$-balls touch, the angle $\phi$ satisfies the equation $\sin\phi=\frac{1}{\lambda r+1}$.  We also know from above that $\phi$ satisfies the equation $\cos\phi=\frac{(R-\frac{1}{\lambda})^2-(R-r-\delta)^2-(r+\frac{1}{\lambda})^2}{2(R-r-\delta)(r+\frac{1}{\lambda})}$.  The curve in figure \ref{circphase} is then derived from the identity $\sin^2\phi+\cos^2\phi=1$.  It is, also, of interest to observe that when the annulus is concentric, i.e. $\delta=.2$, the minimizers for (\ref{l1tv}) are known \cite{yin2007} and coincide with the top line of figure \ref{circphase}.


\section{Square Annulus}

\begin{figure}[h!]
\centering
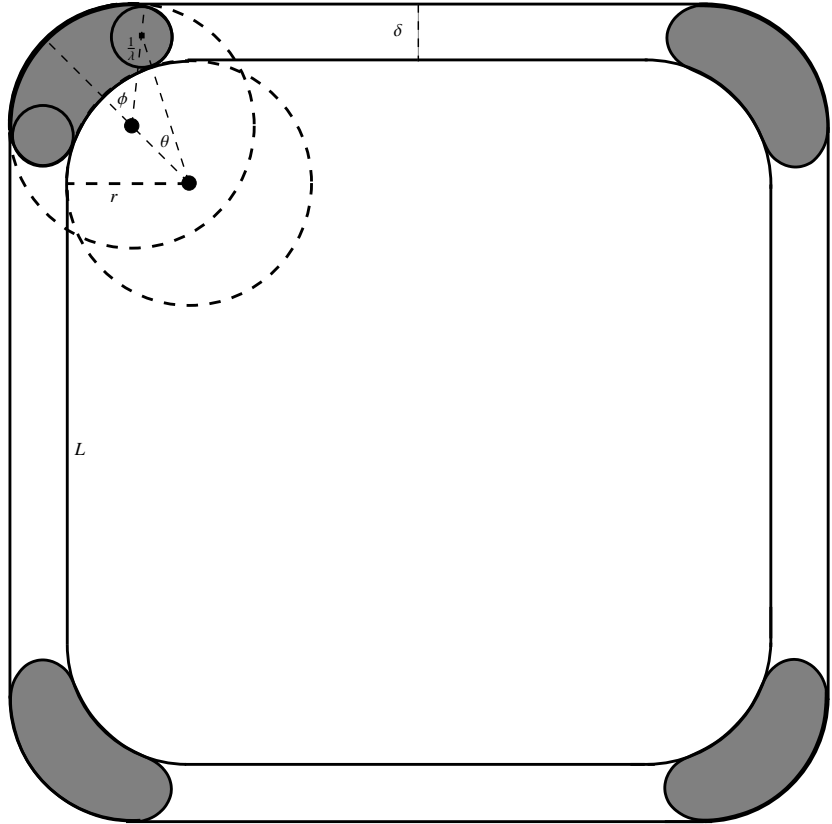
\caption{Square Annulus with angles $\phi$ and $\theta$}\label{squareann}
\end{figure}
For the second example we take $\Omega$ to be the region contained within the `square annulus' shown in figure \ref{squareann}, where the corners are arcs of a circle of radius $r$ with $r>2/\lambda$, $L$ is the length of each side on the inside of the annulus from arc to arc, and $\delta$ is the distance between the straight edges on each side of the annulus.  We then compute (\ref{CE}) for the five choices of $\Sigma_i$:\\
$\Sigma_1=$The large rounded square (the region contained within the outer boundary),\\
$\Sigma_2=\emptyset$,\\
$\Sigma_3=\Omega$,\\
$\Sigma_4=\bigcup\{B_{1/\lambda}(x)\subset\Omega\}$, and\\
$\Sigma_5=$The small rounded square (the region contained within the inner boundary).\\
As long as we note that the curvature at any point of $\partial\Omega$ is less than $\lambda$, we can argue as in the preceding section that these are the only possible minimizers of (\ref{CE}).

This leads to the following equations derived from (\ref{CE}):

\begin{equation}
\label{E6}
E(\Sigma_1,\lambda)=2\pi r+4L+8\delta+\lambda L^2+4\lambda r L+\lambda\pi r^2\end{equation}
\begin{equation}
\label{E7}
E(\Sigma_2,\lambda)=4\lambda\delta L+8\lambda\delta r+4\lambda\delta^2\end{equation}
\begin{equation}
\label{E8}
E(\Sigma_3,\lambda)=4\pi r+8L+8\delta\end{equation}
\begin{equation}
\label{E9}
E(\Sigma_4,\lambda)=8r(\phi+\theta)+\frac{4(\pi-\phi+\theta)}{\lambda}+4\lambda\delta L+8\lambda\delta r+4\lambda\delta^2-4\lambda r^2(\phi-\theta)-4\delta\sqrt{2}(\lambda r+1)\sin\theta
\end{equation}
\begin{equation}
\label{E10}
E(\Sigma_5,\lambda)=2\pi r+4L+\lambda((L+2\delta+2r)^2-4r^2+\pi r^2)
\end{equation}
Where $\phi=\cos^{-1}\left(\frac{2r-\lambda\delta^2}{\sqrt{2}\delta(\lambda r-1)}\right)$ and $\theta=\cos^{-1}\left(\frac{2r+\lambda\delta^2}{\sqrt{2}\delta(\lambda r+1)}\right)$.  The angles $\phi$ and $\theta$ are the angles between the line from the centers of the circles defining the arcs to the corners and the lines from the centers of the circles to the center of one of the eight $\frac{1}{\lambda}$-balls that is tangent to both boundaries, as shown in figure \ref{squareann}.

\begin{figure}[h!]
\centering
\includegraphics[width=\textwidth]{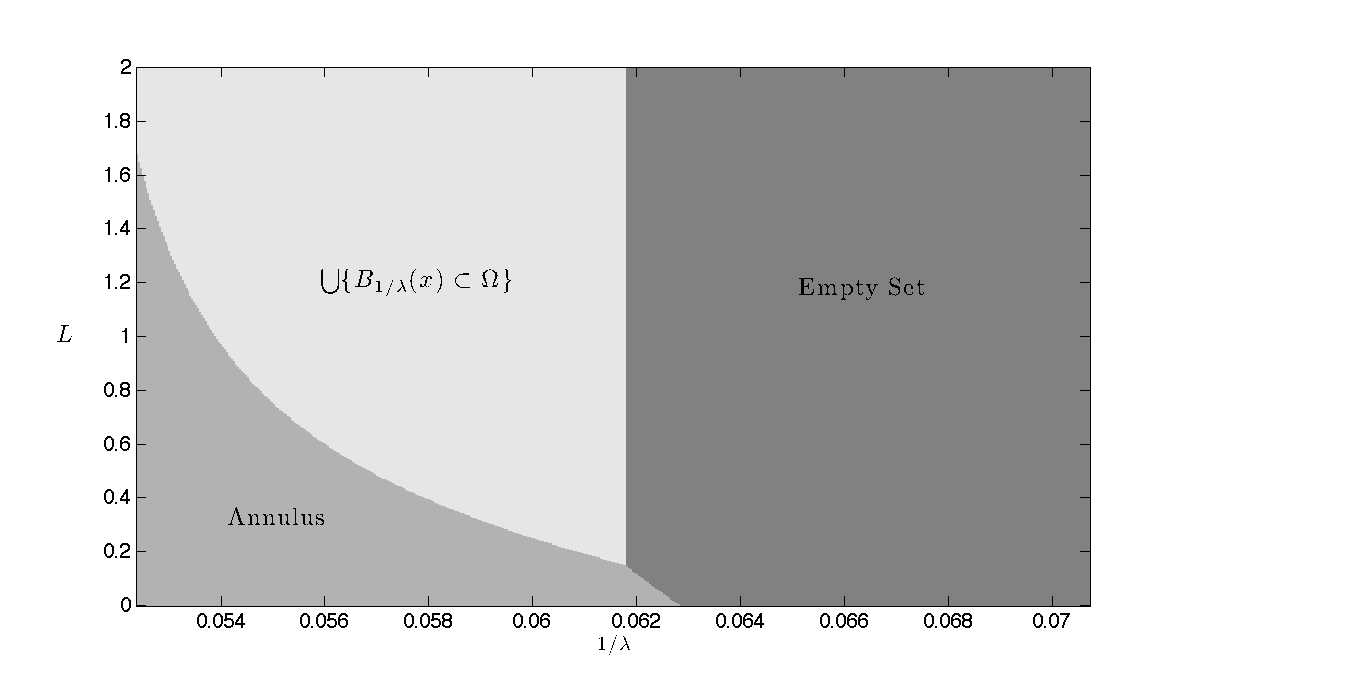}
\caption{This diagram shows the $\Sigma_i$ with minimum value in (\ref{CE}) for the square annulus.}\label{squarephase}
\end{figure}

From equations (\ref{E6})-(\ref{E10}), many examples can be created.  For illustrative purposes, we have chosen $r=1$, $\delta=.1$ and let $L$ and $\lambda$ vary.  We can then compute equations (\ref{E6})-(\ref{E10}) for any values of $L$ and $\lambda$ and determine which equation has minimum value.  Figure \ref{squarephase} summarizes the results for many such choices of $L$ and $\lambda$.  We can then see that in the lower left portion of figure \ref{squarephase} there is a significant region where the entire annulus obtains a lower value for (\ref{CE}) than the union of $\frac{1}{\lambda}$-balls, giving the desired examples.

Note:  We restrict $1/\lambda$ so that it doesn't exceed the widest part of the corners and so that the $1/\lambda$-balls can only be tangent to the two arcs and never to the flat region (in which case equation (\ref{E9}) would be incorrect).  This yields the following bounds $\frac{\delta r+\delta^2}{2r+\delta}<\frac{1}{\lambda}<\frac{\delta}{\sqrt{2}}$.


\section{Dumbbell}

\begin{figure}[h!]
\centering
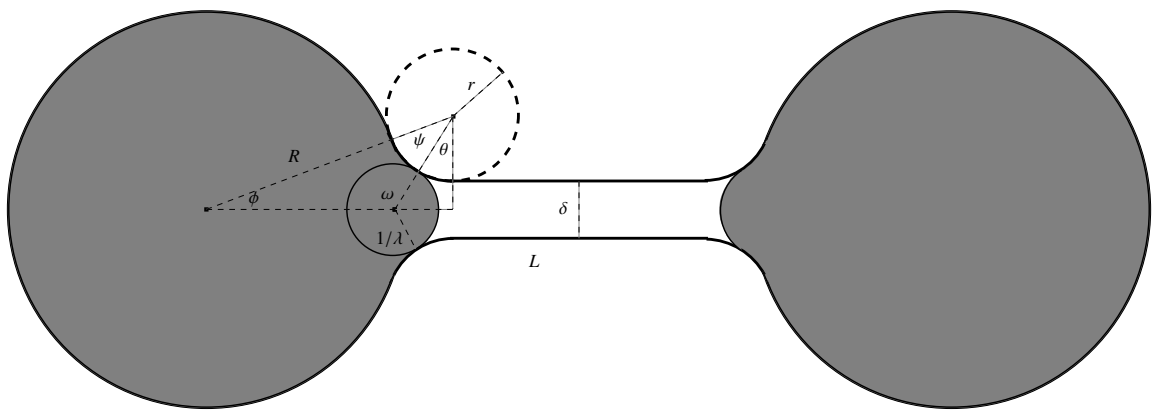
\caption{Dumbbell with angles $\phi$, $\theta$, $\omega$, and $\psi$}\label{dumb}
\end{figure}

For the third example we take $\Omega$ to be the region contained within the `dumbbell' shown in figure \ref{dumb}, where the ends of the dumbbell are circles of radius $R$, the corners between the ends and the `handle' are arcs of a circle of radius $r$ with $2/\lambda<r<R$, $L$ is the length of the `handle' stretching from arc to arc, and $\delta$ is the width of the `handle' with $\delta<2/\lambda$.  We then compute (\ref{CE}) for the three choices of $\Sigma_i$:\\
$\Sigma_1=\emptyset$,\\
$\Sigma_2=\Omega$, and\\
$\Sigma_3=\bigcup\{B_{1/\lambda}(x)\subset\Omega\}$.

Again we can argue that these are the only possible minimizers of (\ref{CE}).  First, it is important to note that we have restricted $\lambda>2/r$, otherwise $\Omega$ would not be a candidate for minimizer because the curvature would exceed $\lambda$.  Since the curvature of $\Omega$ is always smaller than $\lambda$ the minimizer must be contained within $\Omega$ because it is impossible for an arc of $\lambda$-curvature of not more than $\pi$ radians to meet $\Omega$ tangentially from the outside.  One can then argue as before that the only possible nonempty minimizer that is a subset of $\Omega$ is $\bigcup\{B_{1/\lambda}(x)\subset\Omega\}$.

This leads to the following equations:

\begin{equation}
\label{E11}
E(\Sigma_1,\lambda)=2\pi\lambda R^2+\lambda(2r+\delta)\sqrt{(R+r)^2-(r+\delta/2)^2}-2\lambda(\theta+\psi) r^2-2\lambda\phi R^2+\lambda\delta L
\end{equation}
\begin{equation}
\label{E12}
E(\Sigma_2,\lambda)=4(\pi-\phi)R+4(\theta+\psi) r+2L
\end{equation}
\begin{equation}
\label{E13}
\begin{split}
E(\Sigma_3,\lambda)=&4(\pi-\phi)R+4\psi r+\frac{2(\pi-\omega)}{\lambda}+
\lambda\delta L+\lambda(2r+\delta)\sqrt{(R+r)^2-(r+\delta/2)^2}\\&-2\lambda\theta r^2-2(\lambda r+1)(R+r)\sin\psi.
\end{split}
\end{equation}
Where $\phi$, $\theta$, $\omega$, and $\psi$ are all shown in figure \ref{dumb} and are given by $\phi=\sin^{-1}\left(\frac{r+\delta/2}{R+r}\right)$, $\theta=\cos^{-1}\left(\frac{r+\delta/2}{r+1/\lambda}\right)$,  $\omega=\pi/2+\theta$, and $\psi=\pi-\omega-\phi$.

\begin{figure}[h!]
\centering
\includegraphics[width=\textwidth]{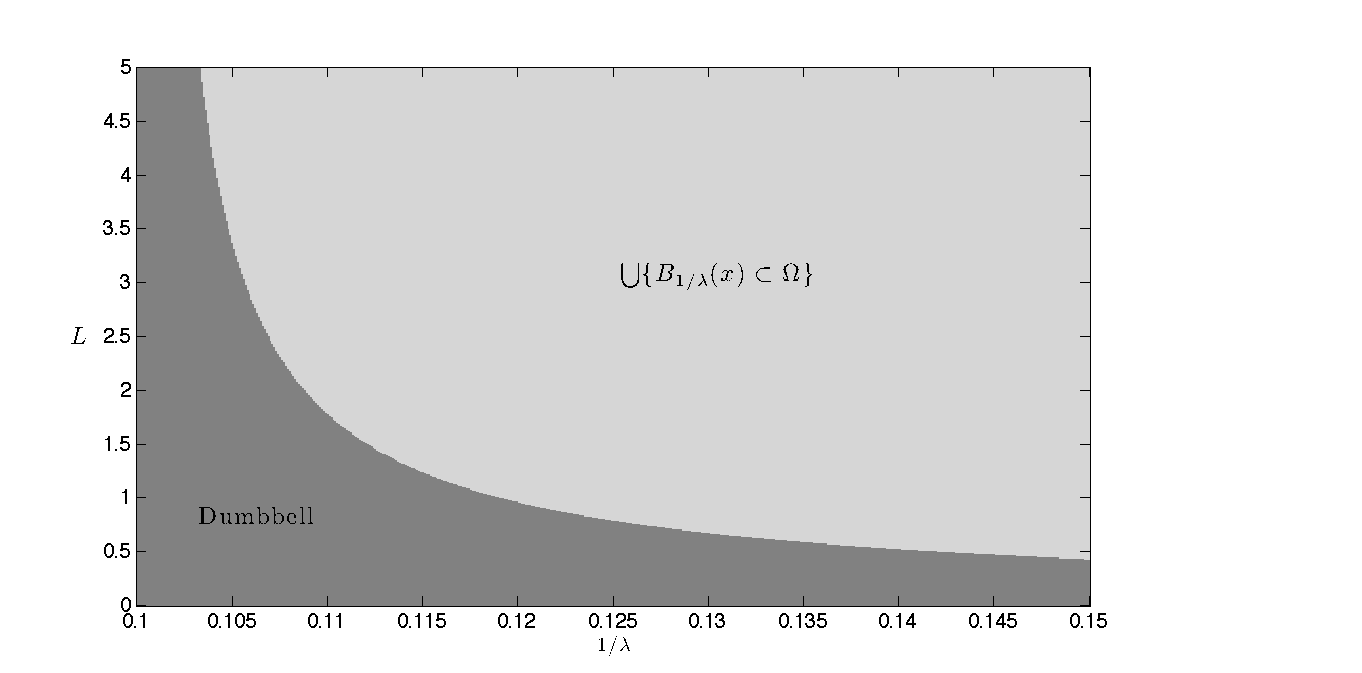}
\caption{This diagram shows the $\Sigma_i$ with minimum value in (\ref{CE}) for the `dumbbell'.}\label{dumbphase}
\end{figure}

From equations (\ref{E11})-(\ref{E13}), many examples can be created.  For illustrative purposes, we have chosen $R=1$, $r=.3$, $\delta=.2$ and let $L$ and $\lambda$ vary.  We can then compute equations (\ref{E11})-(\ref{E13}) for any values of $L$ and $\lambda$ and determine which equation has minimum value.  Figure \ref{dumbphase} summarizes the results for many such choices of $L$ and $\lambda$.  We can then see that in the lower left portion of figure \ref{dumbphase} there is a significant region where the entire `dumbbell' obtains a lower value for (\ref{CE}) than the union of $\frac{1}{\lambda}$-balls, giving the desired examples.

\newpage
\setcounter{page}{1}
\pagestyle{plain}

\bibliographystyle{plain}
\bibliography{l1tv_examples.bib}

\end{document}